

\baselineskip=14pt
\parskip=10pt

\magnification=\magstephalf

\def\1{{\overline{1}}}
\def\2{{\overline{2}}}
\parindent=0pt
\overfullrule=0in

\def\frac#1#2{{#1 \over #2}}
\centerline
{\bf A Quick Empirical Reproof of the  Asymptotic Normality of the  Hirsch Citation Index  }
\centerline
{\bf (First proved by Canfield, Corteel and Savage)}
\bigskip
\centerline
{\it Shalosh B. EKHAD and Doron ZEILBERGER}
\bigskip

{\bf Abstract}: Inspired by Alexander Yong's recent critique of the Hirsch Citation index,
we give an empirical (yet {\bf very convincing!}) reproof of the
{\it asymptotic normality} of the Hirsch citation index (alias size of Durfee square) with respect
to the uniform distribution on the ``sample space'' of integer partitions of $n$. 
This result was first proved rigorously (but with much greater effort!) by the humans 
Rodney Canfield, Sylvie Corteel, and Carla Savage.
In particular, we confirm the
Canfield-Corteel-Savage rigorous evaluation of the average: 
$(\frac{\sqrt{6} \log 2}{\pi})\sqrt{n} +O(1)= 0.5404446 \, \sqrt{n} +O(1)$, and
estimate the variance, numerically, as $0.0811 \, \sqrt{n}+O(1)$, and get estimates 
extremely close to those of the standard Normal Distribution for the first $12$ standardized moments.
We also observe that what Yong calls ``Rodney Canfield's concentration conjecture'',
that asserts that most of the ``mass'' is close to the average, 
follows immediately from the Canfield-Corteel-Savage 1998 result, 
since the {\it variance} is proportional to the average (with a rather small proportionality constant, namely
$0.0811/0.5404446$, that is approximately $0.1501$).
All the results in this article were obtained via straightforward symbol- and number- crunching,
by the aid of a Maple package called {\tt HIRSCH}, that is
available free of charge from {\tt http://www.math.rutgers.edu/\~{}zeilberg/tokhniot/HIRSCH} \quad .

{\bf The Cinderella Story of the Size of the Durfee Square}

Once upon a time there was an esoteric and specialized notion, called ``size of the Durfee square'',
of interest to at most $100$ specialists in the whole world. Then it was kissed by
a prince called Jorge Hirsch ([H]), and became the {\it famous} (and to quite a few people, {\it infamous})
$h$-index, of interest to {\it every} scientist, and scholar, since
it tells you how productive a scientist (or scholar) you are!

When Rodney Canfield, Sylvie Corteel, and Carla Savage wrote their beautiful article, [CCS], proving,
rigorously, by a very deep and intricate analysis, the {\it asymptotic normality} of the random
variable ``size of Durfee square'' defined on integer-partitions of $n$ (as $n \rightarrow \infty$),
with precise asymptotics for the mean and variance, they did not dream that one day
their result should be of interest to everyone who has ever published a paper.

{\bf Alexander Yong's Critique of the h-index}

In the latest issue of the {\it Notices of the American Mathematical Society}, Alexander Yong contributed ([Y])
a very insightful critique of the Hirsch citation index([H]).  Yong mentioned [CCS], but
apparently missed its full significance. In particular, what Yong calls 
``Canfield's concentration conjecture'' is an immediate consequence of the fact, proved in [CCS], that
what now is called the $h$-index is asymptotically normal, and the fact, also proved there,
that the asymptotic variance is proportional to the asymptotic  average.

\vfill\eject

{\bf Why this Redux?}

To be honest, we don't have the patience to follow the intricate analysis of [CCS],
and while we trust them completely, it is nice to find out things by ourselves.
More importantly, we want to describe, via this {\it case-study}, how one can get much
quicker, the same {\it mathematical knowledge}, by combining {\it number-crunching} and
naive {\it symbol-crunching} to get (empirically, but very reliably) limiting distributions
of many families of combinatorial `statistics' (or {\it random-variables}). Often
this method can be used to get fully rigorous results, (see [Z1][Z2]), but with much lesser effort,
one can get empirical proofs. It is also very easy to come up with examples where a fully rigorous
proof is completely beyond the scope of humans, or even computers, and when it is, it is
not worth the efforts!

{\bf What is the h-index (alias Durfee Square)}

Recall that the size of the Durfee square, alias $h$-index, is defined as follows.
For a partition of a positive integer, $n$, $\lambda=(\lambda_1, \dots , \lambda_k)$,
(where $\lambda_1 \geq \dots \geq \lambda_k \geq 1$, and $\lambda_1+ \dots + \lambda_k=n$),
$h(\lambda)$ is the largest $i$ such that $\lambda_i \geq i$.

In  this note we reprove the above-mentioned [CCS] result about the
asymptotic normality of the $h$-index empirically, that immediately implies
the concentration-about-the-mean property,
(what Yong erroneously  thought was only a conjecture, but was in fact a theorem).

{\bf Symbolic Moment Calculus}

In [Z1] (see also [Z2] and [CJZ]) we initiated a symbolic-computational method for
the automatic (and rigorous!) proof of limit laws for many families of combinatorial random variables.
But with much lesser effort, one can always derive the same results {\it empirically} by a combination
of {\it number-crunching} and {\it symbol-crunching} using the very {\it naive} approach that we
will briefly recall.

Let $X_n$ be an infinite sequence of combinatorial families (for example, $\{0,1\}^n$), and let
$f(x)$ be a random variable (for example, the sum of the entries, alias the number of $1$'s).

Define a sequence of polynomials, $C_n(t)$, in a variable $t$, by
$$
C_n(t) := \sum_{x \in X_n} t^{f(x)} \quad ,
$$
(in the above example $C_n(t)=(1+t)^n$),
called the {\it combinatorial generating functions}. Under the {\it uniform distribution}, this
turns into {\it probability generating functions}, by dividing by $C_n(1)$ (alias $|X_n|$)
$$
P_n(t):= \frac{C_n(t)}{C_n(1)} \quad .
$$
(in the above example $C_n(t)=(1+t)^n/2^n$),
To get the expectation, $E_n(f)$, let's call it $a_n$, one simply computes $\frac{d}{dt}P_n(t) \Bigl \vert_{t=1}$
(in the above example $a_n=n/2$). To get the {\it centralized} version, one divides
$P_n(t)$ by $t^{a_n}$, getting
$$
Q_n(t):=\frac{P_n(t)}{t^{a_n}} \quad .
$$
The {\it variance}, $m_2(n)$, is given by
$$
m_2(n)=(t\frac{d}{dt})^2 Q_n(t) \Bigl \vert_{t=1} \quad,
$$
and the higher moments, $k \geq 3$, by
$$
m_k(n)=(t\frac{d}{dt})^k Q_n(t) \Bigl \vert_{t=1} \quad.
$$
Finally, the {\it standardized moments}, $\alpha_k(n)$, are given by
$$
\alpha_k(n)=\frac{m_k(n)}{m_2(n)^{k/2}} \quad .
$$

The random variable $f$ has a {\it limiting distribution} if for every $k>2$,
$$
\beta_k:=\lim_{ n \rightarrow \infty} \alpha_k(n) \quad,
$$
exists. This is usually a {\it continuous} probability distribution, and {\bf very often} 
happens to be the good old
{\it normal distribution} $\frac{1}{\sqrt{2\pi}} e^{-x^2/2}$,
whose moments are $0$ for $k$ odd, and $\frac{k!}{(k/2)!2^{k/2}}$ for $k$ even.
In that case our family of combinatorial random variables is called {\it asymptotically normal}.

In {\it many} cases, this approach can be used to teach the computer to prove {\bf rigorous} results, by getting
either closed-form, or recursive descriptions, of the leading asymptotics, in $n$, of the moments $m_k(n)$,
for {\bf symbolic} $n$ and $k$, and from which one can easily get the leading asymptotics for the standardized moments,
$\alpha_k(n)$. This happens when the $C_n(t)$ are either given explicitly, or via a decent
{\it grand generating function}
$$
F(q,t):=\sum_{n=0}^{\infty} C_n(t) q^n \quad ,
$$
where $F(q,t)$ is more-or-less explicit.

But what if $F(q,t)$ is not so nice? Then we can abandon (alleged) `rigor', and do things {\it empirically}.
Use $F(q,t)$ to crank out the first $10000$ or whatever, polynomials $C_n(t)$, and do all the
above steps up to, say, the $14$-th moment, and estimate asymptotics of $\alpha_k(n)$ for $k \leq 14$ or
whatever. If the leading terms seem to agree with those of the (standard) normal distribution: $1,0,3,0,15,0,105, 0,945, \dots$
we have a very convincing {\it empirical proof} of asymptotic normality. Also, as a bonus we can numerically
estimate asymptotic expressions for the average $a_n$, and the variance $m_2(n)$.

{\bf The Hirsch (formerly Durfee) Polynomials}

Alexander Yong reminds us (Eq. 1 of [A]), about the famous {\bf Euler-Gauss} identity
$$
\prod_{i=1}^{\infty} \frac{1}{1-q^i} = \sum_{k=0}^{\infty} \frac{q^{k^2}}{\prod_{j=1}^{k} (1-q^j)^2} \quad .
$$
This lovely identity was given a pretty combinatorial proof by the Sylvester school (Durfee was
Sylvester's graduate student), that could be found, for example, in George Andrews' partition {\bf bible}([A], pp. 27-28).

That proof immediately implies that the {\it grand generating function}
$$
\sum_{n=0}^{\infty} C_n(t) q^n \quad ,
$$
equals
$$
\sum_{k=0}^{\infty} \frac{q^{k^2} t^k}{\prod_{j=1}^{k} (1-q^j)^2} \quad .
$$
(Note that in order to get the first $N^2$ members of the sequence $C_n(t)$ we only need to take the sum up to $k=N$
and then {\tt taylor} it up to $q^{N^2}$.)

This generating function is not so easy to handle, but Maple can be easily used to crank-out the first
$10000$ terms. The first $6400$ members of the {\it probability generating functions}
(i.e. $C_n(t)/C_n(1)=C_n(t)/p_n$) can be found in {\it Maple input format}, suitable for computer-experimentation,
in the file

{\tt http://www.math.rutgers.edu/\~{}zeilberg/tokhniot/oHIRSCH1} \quad .

The first $100$ terms of the subsequence consisting of perfect squares, i.e. the list of $C_{i^2}(t)$,
for $1 \leq i \leq 100$ can be found  here:

{\tt http://www.math.rutgers.edu/\~{}zeilberg/tokhniot/oHIRSCH2} \quad .

The statistical analysis, as outlined above, can be found in the file

{\tt http://www.math.rutgers.edu/\~{}zeilberg/tokhniot/oHIRSCH3} \quad .

Let's summarize our findings ({\bf Warning}: these are non-rigorous, (but reliable), estimates).
$$
a_n = 0.5404446395\,\sqrt {n}+ 0.085691+ 0.0374788\,{\frac {1}{\sqrt {n}}} + O(\frac{1}{n} ) \quad .
$$
$$
m_2(n)= 0.081057\,\sqrt {n}+ 0.018459- 0.018015\,{\frac {1}{\sqrt {n}}} + O(\frac{1}{n}) \quad .
$$

Finally for the fourth through the eighth standardized even moments, we have:
$$
\alpha_4(n)=3.000000000-.084847493\, \frac{1}{\sqrt{n}} - 0.1071813 \, \frac{1}{n} + O(\frac{1}{n^{3/2}}) \quad ,
$$
$$
\alpha_6(n)= 15.0000000- 12.60947794\,\frac {1}{\sqrt {n}}+ 2.080133651\,\frac{1}{2}+ O(\frac{1}{n^{3/2}}) \quad,
$$
$$
\alpha_8(n)= 105.0000000- 174.8856\, \frac{1}{\sqrt{n}}+104.0909\,\frac{1}{n} + O(\frac{1}{n^{3/2}}) \quad ,
$$
etc. etc. {\bf Very convincing!}.

As already noted in the abstract, the variance is proportional (asymptotically) to
the average. It follows that, since the $h$-index is asymptotically normal, that there is {\bf concentration about the mean}.

{\bf The Maple Package HIRSCH}

Readers are welcome to continue to explore, and generalize, by downloading the
Maple package HIRSCH, already mentioned in the abstract, whose url is:

{\tt http://www.math.rutgers.edu/\~{}zeilberg/tokhniot/HIRSCH} \quad .

{\bf References}

[A] G. Andrews, {\it ``The Theory of Partitions''}, Cambridge University Press, 1998. Originally published
by Addison-Wesley, 1976.

[CCS] E.R. Canfield, S. Corteel and  C.D. Savage,
{\it Durfee polynomials}, Elec. J. Comb. {\bf 5} (1998), R32. (21 pages). Available on-line from: \hfill\break
{\tt  http://www.combinatorics.org/ojs/index.php/eljc/article/view/v5i1r32/pdf } \quad .

[CJZ] E.R. Canfield, S. Janson, and D. Zeilberger, 
{\it The Mahonian Probability Distribution on Words is Asymptotically Normal},
Advances in Applied Mathematics {\bf 46} (2011), 109-124. [Special issue in honor of Dennis Stanton].
Available on-line from: \hfill\break
{\tt http://www.math.rutgers.edu/\~{}zeilberg/mamarim/mamarimhtml/mahon.html} \quad .

[H] J. E. Hirsch, {\it An index to quantify an individual's scientific research output},
Proceedings of the National Academy of Science, {\bf 102 \# 46} (Sept. 1, 2005), 16569-16572.
Available on-line from: \hfill\break
{\tt http://www.pnas.org/content/102/46/16569.full} \quad .

[Y] A. Yong, {\it Critique of Hirsch's Citation Index: a combinatorial Fermi problem},
Notices of the American Mathematical Society {\bf 61 \#9} (Oct. 2014), 1040-1050. Available on-line from: \hfill\break
{\tt www.ams.org/notices/201409/rnoti-p1040.pdf} \quad .

[Z1] D. Zeilberger, {\it The Automatic Central Limit Theorems Generator (and Much More!)},
``Advances in Combinatorial Mathematics: Proceedings of the Waterloo Workshop in Computer Algebra 2008 in honor of Georgy P. Egorychev'', 
chapter 8, pp. 165-174, (I.Kotsireas, E.Zima, eds. Springer Verlag, 2009.) Available on-line from: \hfill\break
{\tt http://www.math.rutgers.edu/\~{}zeilberg/mamarim/mamarimhtml/georgy.html} \quad .

[Z2] D. Zeilberger,
{\it HISTABRUT: A Maple Package for Symbol-Crunching in Probability theory},
The Personal Journal of Shalosh B. Ekhad and Doron Zeilberger,  Aug. 25, 2010. Available on-line from: \hfill\break
{\tt http://www.math.rutgers.edu/\~{}zeilberg/mamarim/mamarimhtml/histabrut.html} \quad .

\vfill\eject

\hrule
\bigskip
Doron Zeilberger, Department of Mathematics, Rutgers University (New Brunswick), Hill Center-Busch Campus, 110 Frelinghuysen
Rd., Piscataway, NJ 08854-8019, USA. \hfill \break
zeilberg at math dot rutgers edu \quad ;  \quad {\tt http://www.math.rutgers.edu/\~{}zeilberg/} \quad .
\bigskip
\hrule
\bigskip
Shalosh B. Ekhad, c/o D. Zeilberger, Department of Mathematics, Rutgers University (New Brunswick), Hill Center-Busch Campus, 110 Frelinghuysen
Rd., Piscataway, NJ 08854-8019, USA.
\bigskip
\hrule

\bigskip
EXCLUSIVELY PUBLISHED IN $\,$: { \tt http://www.math.rutgers.edu/\~{}zeilberg/pj.html} and {\tt arxiv.org}.
\bigskip
\hrule
\bigskip
Oct. 31, 2014

\end